

\input amstex  
\documentstyle{amsppt}
\pagewidth{125mm}
\pageheight{195mm}
\parindent=8mm
\frenchspacing\tenpoint
\NoRunningHeads

\batchmode
  \font\gothic=eufm10
  \font\footgothic=eufm8
\errorstopmode

\ifx\gothic\nullfont
  \define\uu{\operatorname{u}}
\else
  \define\uu{\operatorname{\hbox{\gothic u}}}
\fi
\ifx\footgothic\nullfont
  \define\footuu{\operatorname{u}}
\else
  \define\footuu{\operatorname{\hbox{\footgothic u}}}
\fi

\hyphenation{Ca-la-bi}

\def\pretitle{\noindent\hskip23pc\vbox to0pt{%
                    \vskip-31pt\advance\hsize-23pc\parindent0pt
     DUK-M-94-06    \hfil\break
     November, 1994 \par\vss}\vskip3pc}

\topmatter
\title Mirror Symmetry and Moduli Spaces\\
of Superconformal Field Theories \endtitle
\author David R. Morrison \endauthor
\address
Department of Mathematics,
Duke University,
Durham, NC 27708-0320 USA
\endaddress
\email drm\@math.duke.edu \endemail
\thanks Research partially supported by  National Science Foundation
Grants DMS-9103827, DMS-9304580 and DMS-9401447, and by
an American Mathematical Society Centennial Fellowship. \endthanks
\endtopmatter

\define\C{{\Bbb C}}
\define\R{{\Bbb R}}
\define\Z{{\Bbb Z}}

\define\CP{{\Bbb C \Bbb P}}

\define\Aut{\operatorname{Aut}}
\define\Diff{\operatorname{Diff}}
\define\GL{\operatorname{GL}}
\define\Hom{\operatorname{Hom}}
\define\Id{\operatorname{Id}}
\define\SU{\operatorname{SU}}
\define\Sym{\operatorname{Sym}}
\define\U{\operatorname{U}}

\define\eqref#1{(#1)}


\define\bosonicaction{1}
\define\Msigma{2}
\define\decomposition{3}
\define\topological{4}
\define\pathint{5}
\define\pathtermbis{6}
\define\asympt{7}
\define\jprime{8}
\define\MNtwo{9}
\define\cKcone{10}
\define\mirrormap{11}


\define\AGM{1}
\define\small{2}
\define\AMone{3} 
\define\AMtwo{4} 
\define\AMthree{5} 
\define\BE{6}
\define\Besse{7}
\define\Bogomolov{8}
\define\Calabi{9}
\define\CdGP{10}
\define\cls{11}
\define\DS{12} 
\define\dixon{13}
\define\voa{14}
\define\GSW{15}
\define\gp{16}
\define\Gromov{17}
\define\GS{18}
\define\lvw{19}
\define\MS{20}
\define\guide{21}
\define\compact{22}
\define\beyond{23}
\define\predictions{24}
\define\Narain{25}
\define\NSW{26}
\define\RT{27}
\define\Tian{28}
\define\Todorov{29}
\define\Viehweg{30}
\define\WBerkeley{31}
\define\topsigma{32}
\define\topgrav{33}
\define\Wmirror{34}
\define\phases{35}
\define\YauICM{36}

\document

{\it Mirror symmetry}\/ is the remarkable discovery in string theory that
certain ``mirror pairs'' of Calabi--Yau manifolds
apparently
produce isomorphic physical theories---related by an isomorphism
which reverses the sign of a certain
quantum number---when used as backgrounds for string
propagation \cite{\dixon, \lvw, \cls, \gp}.
The sign reversal in the isomorphism has profound effects
on the geometric interpretation of the pair of physical theories.
This leads
to startling predictions that certain geometric invariants of
one Calabi--Yau
manifold (essentially the numbers of holomorphic $2$-spheres of
various degrees) should be
related to a completely different set of geometric invariants of
the mirror partner (``period'' integrals of holomorphic forms).

We will discuss the applications of this mirror symmetry principle to the
study of the moduli spaces
of two-dimensional conformal field
theories with $N{=}(2,2)$ supersymmetry.  Such theories depend on finitely
many parameters, and for a large class of these theories the
parameters admit a clear geometric interpretation.  To circumvent the
difficulties of trying to treat path integrals in a mathematically
rigorous manner,
we shall simply {\it define}\/ the moduli spaces in terms of these geometric
parameters.  Other interesting physical quantities---the
``topological''
correlation functions---can then also be defined
as asymptotic series whose coefficients have geometric meaning.  The
precise forms of the
definitions are motivated by path integral arguments.

Mirror symmetry predicts some unexpected identifications between these
moduli spaces, and serves as a powerful tool for understanding their
structure.  Perhaps the most striking consequence is the prediction
that the moduli spaces can be analytically continued beyond the original
domain of definition, into new regions, some of
 which parameterize conformal field
theories that are related not to the original Calabi--Yau
manifold, but rather to close cousins of it which differ by  simple
topological transformations.

In preparing this report, I have drawn on a considerable body of earlier work
\cite{\AGM--\AMthree, \guide--\predictions},
much of
which was
 collaborative.  I would like to thank my colleagues and
collaborators
Paul Aspinwall, Robert Bryant, Brian Greene, Sheldon Katz,
Ronen Plesser, and Edward Witten for their contributions.

\head 1. The physics of nonlinear $\sigma$-models \endhead

We begin by describing nonlinear $\sigma$-models from the
point of view of physics (see \cite{\GSW} and the references therein),
and giving a geometric interpretation to the parameters
which appear in the theory.
The starting data for constructing a
nonlinear $\sigma$-model consists of a
compact manifold $X$,
a Riemannian metric $g_{ij}$  on $X$, and a class
$B\in H^2(X,\R/\Z)$
(which we represent as a closed, $\R/\Z$-valued $2$-form,
i.e., a collection of closed, $\R$-valued $2$-forms on the sets of an open
cover of $X$ which differ by $\Z$-valued forms on overlaps).
The bosonic version of the nonlinear $\sigma$-model is then specified,
in the Lagrangian formulation,
by the $\C/\Z$-valued (Euclidean) action
which assigns to each sufficiently smooth
 map $\phi$
from an oriented Riemannian $2$-manifold $\Sigma$ to $X$
the quantity\footnote{We suppress
the string coupling constant, and use a normalization in which
the action appears as $\exp(2\pi i{\Cal S})$
in the path integrals for correlation functions.}
$${\Cal S}[\phi]:=
i\int_\Sigma \|d\phi\|^2\,d\mu+\int_\Sigma \phi^*(B),
\tag\bosonicaction$$
where the norm $\|d\phi\|$ of $d\phi\in\Hom(T_\Sigma,\phi^*(T_X))$
is determined
from the Riemannian metrics on $X$ and on $\Sigma$.

There is a variant of this theory in which additional
fermionic terms are added to \eqref{\bosonicaction} to produce
an action which is invariant under at least one supersymmetry
transformation. (We will not write
the fermionic terms in the action explicitly, as they do not enter
into our analysis of the parameters.)
The supersymmetric form of the
action is also  invariant under {\it additional}\/ supersymmetry
transformations when the geometry is
restricted in certain ways---if the metric
is K\"ahler then the theory has what is called
$N{=}(2,2)$ supersymmetry, while if the metric is hyper-K\"ahler
then the supersymmetry algebra is extended to $N{=}(4,4)$.

A nonlinear $\sigma$-model describes a consistent background for string
propagation only if it is conformally invariant.
The possible failure of conformal invariance is measured by the so-called
``$\beta$-function'' of the theory, and a perturbative calculation
yields the result that the one-loop contribution to this $\beta$-function
is proportional to the Ricci tensor of the  metric.
This makes Ricci-flat metrics---those with vanishing Ricci tensor---into
good candidates for producing conformally
invariant $\sigma$-models.
In fact, supersymmetric $\sigma$-models whose Ricci-flat metric
is in addition hyper-K\"ahler
are believed to be conformally
invariant, as are bosonic $\sigma$-models whose metric is flat.

When the supersymmetry algebra of the theory cannot be extended
as far as $N{=}(4,4)$,
the Ricci-flat theories fail to be conformally invariant.
However, when the Ricci-flat metrics are K\"ahler (i.e., when the
theory has $N{=}(2,2)$ supersymmetry), we can deduce some of the
properties of the conformally invariant theory by a careful study
of the Ricci-flat theories.  This works as follows:  renormalization
produces a flow on the space of metrics, and along a trajectory which
begins at a Ricci-flat K\"ahler metric, the metric is expected to remain
K\"ahler with respect to a fixed complex structure on $X$, and the K\"ahler
class of the metric is not expected to change.
Thus, if there is
a conformally invariant theory in the same universality class as
this trajectory, i.e.,
if there is
a fixed point of the flow which lies in the trajectory's closure,
then any property
of the conformal theory which depends only on the complex structure,
the  K\"ahler class,
and the $2$-form $B$ can be calculated anywhere along the trajectory,
including the initial, Ricci-flat theory.  Furthermore, every
 Ricci-flat K\"ahler metric
whose K\"ahler class is sufficiently deep within the K\"ahler cone
 is expected to
determine a unique conformally invariant theory (which lies in
the same universality class).

We can thus define a first approximation to the parameter space for
$N{=}(2,2)$ superconformal field theories as follows (cf.\ \cite{\compact}).
Fix a compact manifold $X$, and define
the {\it  one-loop semiclassical nonlinear $\sigma$-model moduli
space}\/ of $X$ to be
$${\Cal M}_\sigma:=\{(g_{ij},B)\}/\Diff(X),
\tag{\Msigma}$$
where
 $g_{ij}$ runs over the set of
Ricci-flat metrics which are K\"ahler for some complex structure on $X$,
 $B$ is an element of $H^2(X,\R/\Z)$, and
 $\Diff(X)$ denotes the diffeomorphism
group of $X$.
Manifolds for which ${\Cal M}_\sigma$ is nonempty (that is, those which
admit a Ricci-flat K\"ahler metric) are called {\it Calabi--Yau manifolds}.
The $2$-form $B$ should be regarded as some sort of ``extra structure''
(cf.\ \cite{\guide})  which supplements the choice of metric.

It is important to keep in mind that the space ${\Cal M}_\sigma$
 is only an  approximation to the moduli space of conformal field
theories, for several reasons:
\item{$\bullet$} As already mentioned,
not every pair $(g_{ij},B)$ is expected to determine
a conformal field theory, only those
whose K\"ahler class is sufficiently deep within the K\"ahler cone.
\item{$\bullet$} There may be analytic continuations of the space of
conformal field
theories beyond the domain where the theories have  a $\sigma$-model
interpretation. (We will see this in more detail in section 6.)
\item{$\bullet$} There may be points of ${\Cal M}_\sigma$ which define
isomorphic conformal field theories, even though they do not define
isomorphic $\sigma$-models.  This phenomenon was first observed in the
case in which $X$ is a torus of real dimension
 $2d$, and $g_{ij}$ is a flat metric
\cite{\Narain, \NSW}:  in this case,
${\Cal M}_\sigma=\Gamma_0\backslash{\Cal D}$, where ${\Cal D}$ is a
certain symmetric space and $\Gamma_0=\Lambda^2\Z^{2d}\rtimes\GL(2d,\Z)$,
while the actual
moduli space of conformal field theories
takes the form $\Gamma\backslash{\Cal D}$
for some $\Gamma$ containing the integral orthogonal group
$\operatorname{O}(\Z^{2d,2d})$ (in which $\Gamma_0$ is a parabolic subgroup).

\noindent
In spite of these limitations, ${\Cal M}_\sigma$
 provides a good arena for
formulating a mathematical version of the theory, based on definitions
using asymptotic expansions.

\head 2. The  correlation functions \endhead

The  correlation functions of these quantum field theories
will depend on the parameters in the action functional.  If we construct
a vector bundle over the moduli space whose fiber over a particular point
is the
Hilbert space of operators in the theory labeled by that point,
then the correlation functions can be regarded as multilinear maps
from this bundle to the complex numbers.
  These maps and their dependence on parameters can
be studied by means of a semiclassical analysis, at least in a certain
``topological'' sector of the theory.
(In this sector, the dependence of the correlation functions on the metric
will always be a dependence on the K\"ahler class alone.)

The semiclassical properties of the $N{=}(2,2)$ theory
are calculated in terms of the set of
stationary values for the action
\eqref{\bosonicaction}.
To find these, we
pick a complex structure on $\Sigma$ which makes
its Riemannian metric K\"ahler, and which is compatible with its
orientation.  Then the first term in the action \eqref{\bosonicaction} can be
rewritten using the formula:
$$\int_\Sigma \|d\phi\|^2\,d\mu=\int_\Sigma \|\bar\partial\phi\|^2\,d\mu
+\int_\Sigma \phi^*(\omega),
\tag\decomposition$$
where
$\bar\partial\phi\in\Hom(T_\Sigma^{(1,0)},\phi^*(T_X^{(0,1)}))$
is determined by the complex structures, and where
 $\omega$ is the K\"ahler form of the metric $g_{ij}$ on $X$.
{}From this formula  it is clear that the stationary values are the holomorphic
maps, i.e., those with $\bar\partial\phi\equiv0$.
Furthermore, the action (\bosonicaction) evaluated on such a stationary value
is the
quantity
$$i\int_\Sigma \phi^*(\omega)+\int_\Sigma \phi^*(B)\in\C/\Z,
\tag\topological$$
which depends only on the homology class $\eta$ of the map $\phi$.

The path integral describing this quantum field theory has bosonic part
$$\int {\Cal D}\phi\,e^{2\pi i\,{\Cal S}[\phi]},
\tag{\pathint}$$
and the correlation functions are calculated by inserting operators into
this expression (see for example Witten's address at the Berkeley ICM
\cite{\WBerkeley}).
Such path integrals are of course problematic for mathematicians,
but it is possible to use the outcome of the path integral manipulations as
a basis for mathematical definitions.

To analyze these  correlation functions,
 we break the path integral into a sum over homology
classes.  This produces an asymptotic expansion which is expected to
converge for metrics whose K\"ahler class is sufficiently deep within
the K\"ahler cone.
The terms in the asymptotic expansion
are themselves path integrals whose bosonic parts are the integrals of
$\exp(2\pi i\int_\Sigma \|\bar\partial\phi\|^2\,d\mu)$ over all maps
of class $\eta$
(with operators inserted),
weighted by the exponential of $2\pi i$ times the classical action
\eqref{\topological}.
For certain of the correlation functions,
these ``coefficient'' path integrals
can in turn be evaluated by
the methods of topological field theory (cf.\ \cite{\topsigma, \Wmirror}):
upon modifying the fermionic terms in the action and introducing
a parameter $t$, the path integral with bosonic part
$$\int_{[\phi]=\eta}{\Cal D}\phi\,e^{2\pi i\,t
\int_\Sigma \|\bar\partial\phi\|^2\,d\mu}
\tag\pathtermbis$$
and ``topological'' operator insertions
becomes independent of $t$.  This integral can then be evaluated by the
method of stationary phase, which reduces it to a finite-dimensional
integral over the set of stationary maps in class $\eta$.
Rigorous mathematical definitions for such ``topological''
correlation functions can be based
on these finite-dimensional integrals, following ideas of
Gromov \cite{\Gromov} and Witten \cite{\topsigma, \topgrav}.
See \cite{\MS}, \cite{\RT}, or Kontsevich's address at this Congress for an
account of these definitions and their properties.

In short, the physical quantities which can be calculated (by physicists) or
defined (by mathematicians) using
topological field theory
will take the general form
$$\sum_{\eta\in H_2(X,\Z)} c_\eta\,
e^{2\pi i\,\langle B+i\omega,\eta\rangle}.
\tag\asympt$$
Notice that the only dependence on the metric is through the complex structure
and the K\"ahler class
$\omega$.
The coefficient $c_\eta$ will depend on the set of all holomorphic maps
in class $\eta$, and may well depend on the complex structure of $X$.
(It also depends on the behavior of the fermionic terms in the
action which we have suppressed.)
The key property of interest here is the holomorphic dependence of
these functions on parameters:  the coefficients $c_\eta$ depend
holomorphically on the complex structure, and the dependence
of \eqref{\asympt} on
$B+i\omega$ is also holomorphic (provided that the series converges
and that $H^{2,0}(X_{\Cal J})=\{0\}$.)

\head 3. Mirror symmetry \endhead

The analysis of the previous sections ultimately derives from the
specific form of our physical theory, which is based on the geometry
of Ricci-flat metrics on $X$.  We now adopt a somewhat more abstract
point of view, and consider the structure of $N{=}(2,2)$ superconformal
field theories {\it per se}.

The algebraic approach to conformal field theories---which treats them
as unitary representations of the Virasoro algebra---has been
extensively studied in the mathematics literature (cf.~\cite{\voa}, for
example).
When the theories
are supersymmetric, the algebra which acts on the representation
can be enlarged.  The enlargement relevant here is the
$N{=}2$ superconformal algebra
 (for which a convenient reference is \cite{\lvw}).
This is a super extension of the Virasoro algebra whose even part
contains a $\uu(1)$-subalgebra in addition to the Virasoro algebra itself.
{}From this algebraic point of view, an $N{=}(2,2)$ superconformal
field theory is simply
a unitary representation of two commuting copies of this algebra;
there is thus
an induced representation of the subalgebra $\uu(1)\times\uu(1)$.

The deformations of these representations have been analyzed in the physics
literature \cite{\DS, \dixon}.  The infinitesimal deformations can
be identified with the finite-dimensional kernel $V$ of a certain operator,
and it is argued in \cite{\DS, \dixon}
that there should be no obstructions to deforming
in the directions corresponding to $V$.\footnote{The arguments in
\cite{\DS} and \cite{\dixon}
involve more of the physical structure than is present in the purely
algebraic formulation we are discussing here.  It would be desirable to
have a purely algebraic proof of this statement.}

The $\uu(1)\times\uu(1)$ manifests itself on $V$ in the following
way:  there are two commuting
complex structures ${\Cal J}$ and ${\Cal J}'$ on
$V$, each of which determines a natural representation of
 $\uu(1)$ on $V\otimes\C$
(with respect to which half of the charges\footnote{For a representation
$\rho$ of $\footuu(1)\cong i\R$, the eigenvalues of $\rho(i)$ are called the
{\it charges}\/ of the representation.} are $+1$ and half are $-1$).
The two complex structures together determine a representation of
 $\uu(1)\times\uu(1)$  on $V\otimes\C$,
and we can decompose $V\otimes\C$ into four complex subspaces
$V^{\pm1,\pm1}$ according to the $\uu(1)\times\uu(1)$ charges.

If we use ${\Cal J}$ to put a complex structure on $V$ and call the
resulting space $V_{\Cal J}$, then we can write
$V_{\Cal J}=V^{1,1}\oplus V^{1,-1}$.  From this point
of view, since ${\Cal J}'$ has eigenvalues $\pm i$, respectively, on
the two summands while ${\Cal J}$
is simply multiplication by $i$,
 we can identify the summands as
$V^{1,1}=\ker({\Cal J}{\Cal J}'-\Id)\subset V_{\Cal J}$  and
$V^{1,-1}=\ker({\Cal J}{\Cal J}'+\Id)\subset V_{\Cal J}$.

A {\it mirror isomorphism}\/ between two $N{=}(2,2)$ superconformal
field theories is an isomorphism which reverses the sign of {\it one}\/ of the
$\uu(1)$ charges.  If it is the second $\uu(1)$ charge
which is reversed, then the isomorphism will map
$V^{\pm1,\pm1}$ to $V^{\pm1,\mp1}$,
and will interchange the factors in the decomposition
$V_{\Cal J}=V^{1,1}\oplus V^{1,-1}$.

A mirror isomorphism must preserve {\it all}\/ correlation functions,
not just the topological ones.  It particular, it preserves the bilinear
form on $V$ which corresponds to the so-called Zamolodchikov metric
on ${\Cal M}_\sigma$.  Thanks to the preservation of this metric,
a mirror isomorphism at a single point can always be extended to a local
isometry between the moduli spaces.
There will also be a compatible
isomorphism of the bundles of Hilbert spaces which maps the topological
correlation functions from one theory to those of the other, but
because of the sign change in the $\uu(1)$ charge, the geometric
interpretations of these correlation functions may be rather different.
For example, Candelas, de la Ossa, Green, and Parkes
\cite{\CdGP}
used a mirror isomorphism to assert that a correlation function which
they could compute exactly (using period integrals) as
$$5+2875\,{q\over1-q}+609250\,{2^3q^2\over1-q^2}+
317206375\,{3^3q^3\over1-q^3}+
242467530000\,{4^3q^4\over1-q^4}+
\cdots$$
should coincide with a generating function
of the form \eqref{\asympt} in which the coefficients represent the
numbers of holomorphic $2$-spheres of various degrees on a quintic hypersurface
in $\CP^4$.
  (See \cite{\guide} or the Givental's address at this Congress
for some of the mathematical aspects of this generating function.)

\head 4. Local analysis of the $\sigma$-model moduli space \endhead

The abstract description of the deformations of $N{=}(2,2)$ theories
can be made very concrete for $\sigma$-models, where it reveals
the local structure of the space ${\Cal M}_\sigma$.
The set of first-order variations $\delta g$
of a fixed Riemannian metric $g_{ij}$ on $X$ can be
identified with the space of symmetric contravariant $2$-tensors
$\Gamma(\Sym^2T_X^*)$.  If $X$ is compact and $g_{ij}$ is Ricci-flat,
then according to a theorem of Berger and Ebin
 \cite{\BE}  the space of first-order variations\footnote{It follows from
the theorem of Bogomolov \cite{\Bogomolov}, Tian \cite{\Tian} and
Todorov \cite{\Todorov} that first-order variations
can always be extended to deformations of the metric.} $\delta g$
(modulo $\Diff(X)$)
which preserve the Ricci-flat condition
can be identified
with the kernel of the Lichnerowicz Laplacian $\Delta_L$
acting on $\Gamma(\Sym^2T_X^*)$.
On the other hand, the set of first-order variations
 $\delta B$ of the $2$-form
$B$ can be identified with the space of harmonic $2$-forms
$\ker\Delta\subset\Gamma(\Lambda^2T_X^*)$.  Since the Lichnerowicz Laplacian
on $2$-forms coincides with the ordinary Laplacian, the combined
 contravariant $2$-tensor $\delta g+\delta B\in \Gamma(\bigotimes^2T_X^*)$
satisfies $\Delta_L(\delta g+\delta B)=0$.  We can thus identify
 the tangent space to ${\Cal M}_\sigma$ at $(g_{ij},B)$
with
$\ker \Delta_L\subset\Gamma(\bigotimes^2T_X^*)$.

Let us assume that the holonomy of $g_{ij}$ takes its ``generic'' value
for Ricci-flat K\"ahler metrics, namely $\SU(n)$, $n\ge3$
(where $n:=\dim_{\R}X$).  In this case,
the two complex structures which we are expecting from our abstract
analysis can be described as follows.
First,
if we fix a complex structure\footnote{When the holonomy is $\SU(n)$,
$n\ge3$, there
are precisely two such complex structures:  ${\Cal J}$ and $-{\Cal J}$.}
${\Cal J}$  on $X$ with respect to which $g_{ij}$
is K\"ahler, there is an induced
operator ${\Cal J}$ on $\Gamma(\bigotimes^2T_X^*)$
defined
by
${\Cal J}h(x,y):=h(x,{\Cal J}y)$.
This new operator ${\Cal J}$ commutes with $\Delta_L$, and so induces an
operator on the tangent
space $\ker{\Delta_L}$ of ${\Cal M}_\sigma$ whose square is $-\Id$,
that is, a complex structure on $\ker\Delta_L$.

The second complex structure ${\Cal J}'$ on $\ker\Delta_L$ is much less
obvious.
It can be characterized by the property that the product
${\Cal J}{\Cal J}'$ acts as $-\Id$ on the space of symmetric,
skew-Hermitian tensors, and as $+\Id$ on the space of tensors which
are either Hermitian or skew-symmetric.  Explicitly, ${\Cal J}'$
 can be defined by the formula
$${\Cal J}'h(x,y):=
{1\over2}\left(-h(x,{\Cal J}y)+h(y,{\Cal J}x)
+h({\Cal J}x,y)+h({\Cal J}y,x)\right).
\tag{\jprime}$$

Using ${\Cal J}$ to put a complex structure on $\ker\Delta_L$,
we can identify $V^{1,-1}$ with
the space of symmetric, skew-Hermitian tensors in
$\ker\Delta_L$.
This space corresponds to that part of the moduli space
of metrics which is obtained by varying the complex structure
(cf.\ \cite{\Besse, Chapter 12}).  The operator ${\Cal J}$ preserves that
space, and induces the usual complex structure on it.
In fact, under our assumptions about the holonomy,
 the complex structure can be
varied freely and we have
 $V^{1,-1}\cong H^1(T_X^{(1,0)})$,
the latter being the space of first-order variations of complex structure.

We can similarly identify $V^{1,1}$ as the
space consisting of
tensors which are either Hermitian or skew-symmetric; on this space,
the operator ${\Cal J}$ mixes symmetric and skew-symmetric forms,
so does not have a classical interpretation in terms of metrics alone.
The parameters associated to this part of the deformation space are
of the form $B+i\omega$, and $V^{1,1}\cong H^{1,1}(X_{\Cal J})\cong H^2(X,\C)$
(under our assumption that the holonomy is $\SU(n)$, $n\ge3$).

A mirror isomorphism between Calabi--Yau manifolds $X$ and $Y$ thus
identifies the space of complex deformations of $X$ with the space
of complexified K\"ahler deformations of $Y$, and vice versa
(at least when the holonomy is ``generic'').

\head 5. Global analysis of the $\sigma$-model moduli space \endhead

The moduli space of Ricci-flat
metrics (and hence the
nonlinear $\sigma$-model moduli space) can be analyzed globally as
well as locally.  To carry this out,
we introduce a related space which includes
a choice of complex structure.
Define
$${\Cal M}_{N{=}2}:=\{(g_{ij},B,{\Cal J})\}/\Diff(X)
\tag{\MNtwo}$$
where ${\Cal J}$ ranges over the complex structures on $X$ with
respect to which $g_{ij}$ is K\"ahler.
The holonomy group of the metric $g_{ij}$ is necessarily contained
in the $\SU(n)$ specified by ${\Cal J}$.
The fibers of the natural map
${\Cal M}_{N{=}2}\to{\Cal M}_\sigma$
depend on this holonomy group, and can be described
as the set of $\U(n)$'s which lie
between the holonomy group and $\operatorname{O}(2n)$. Some examples:
\item{1.} If the holonomy is $\SU(n)$, $n{\ge}3$, then the
fiber consists of two points. (This is the ``generic'' case.)
\item{2.} If the holonomy is $\operatorname{Sp}(n/2,\C)$, then the fiber is
$\CP^1$. (This is the case of an indecomposable hyper-K\"ahler manifold,
such as a K3 surface.\footnote{K3 surfaces are ``self-mirror,'' and the
mirror map induces an automorphism of ${\Cal M}_\sigma$.
Thus, as in the case of a torus, the moduli space of
 conformal field theories of this type is a nontrivial quotient of
${\Cal M}_\sigma$ (cf.\ \cite{\AMthree}, where this quotient is
determined precisely).})

\noindent
The real dimension of the fiber is always $\dim_{\R}H^{2,0}(X_{\Cal J})$.

The structure of the space ${\Cal M}_{\text{N{=}2}}$ can be determined
from the natural map
${\Cal M}_{N{=}2}\to{\Cal M}_{\text{complex}}:=\{{\Cal J}\}/\Diff(X)$.
By the theorems of Calabi \cite{\Calabi} and Yau \cite{\YauICM}, the
fibers of this map
take the form ${\Cal K}_{\C}(X_{\Cal J})/\Aut(X_{\Cal J})$,
where ${\Cal K}_{\C}(X_{\Cal J})$ is
 the {\it complexified K\"ahler cone}\/\footnote{This definition
differs slightly from ones we have given elsewhere \cite{\compact, \beyond}.}
$$
{\Cal K}_{\C}(X_{\Cal J}):=
\{B+i\omega\in H^2(X,\C/\Z) \ |\ \omega\in{\Cal K}_{\Cal J}\},
\tag{\cKcone}$$
${\Cal K}_{\Cal J}$ being the set of K\"ahler classes on
$X_{\Cal J}$, and $\Aut(X_{\Cal J})$ being the group of
holomorphic automorphisms.
It is this fact which gives us access to global information about the
conformal field theory
moduli space, since the moduli space of complex structures  can
be studied by the methods of algebraic geometry.  For example,
by a theorem of Viehweg \cite{\Viehweg} the
subspace ${\Cal M}_{\text{complex}}^{\Cal L} \subset
{\Cal M}_{\text{complex}}$ consisting of all complex structures
polarized with respect to a fixed class
 ${\Cal L}$ is a quasi-projective variety, i.e.,
the complement of a finite number of compact subvarieties in a compact
complex manifold.  (And the spaces ${\Cal M}_{\text{complex}}^{\Cal L}$
are open subsets of ${\Cal M}_{\text{complex}}$
when $H^{2,0}(X_{\Cal J})=\{0\}$.)
In contrast, although ${\Cal K}_{\C}$ has a canonical complex structure when
$H^{2,0}(X_{\Cal J})=\{0\}$, it is typically a rather small domain.

Note that the expected condition for a given pair $(g_{ij},B)$
to determine a conformal field theory was stated in terms of the
K\"ahler class only and was valid for every choice of
 complex structure.
Thus, the global description of the complex structures should be
valid for the conformal field theory moduli space itself.
On the other hand, the complexified K\"ahler directions are
subject to modification.

\head 6. Beyond the K\"ahler cone \endhead

We now apply the mirror symmetry principle to study the
moduli space in the case in which
 the holonomy of the Ricci-flat
metrics on $X$ is
$\SU(n)$, $n\ge3$.

Suppose that a mirror partner $Y$ is known for $X$.  The mirror map
between the moduli spaces ${\Cal M}_\sigma(X)$ and ${\Cal M}_\sigma(Y)$
will certainly be well-defined at points corresponding to metrics whose
K\"ahler class is
 sufficiently
deep within the K\"ahler cone, but in general we can only expect a
 partially defined, local isomorphism between these spaces.  However,
because of the global nature of the complex structure space
${\Cal M}_{\text{complex}}(Y)$, we can deduce the structure of
the K\"ahler moduli space ${\Cal K}_{\C}(X)$ from even a {\it local}\/
knowledge
of the mirror map.  In principle, the mirror map should be determined
essentially uniquely from the structure of the Zamolodchikov metric,
once the derivative of the map is known at a single point.
In practice, it is easier to approach the construction of the
mirror map in other ways (based on the topological correlation
functions) which determine it up to a finite number
of unknown parameters.  Even those parameters can often be
determined.  (See \cite{\predictions} for a recent review of this
problem.)

This comparison of structure between K\"ahler and complex moduli spaces has
been
carried out in \cite{\AGM, \small}
for  cases in which a mirror partner is known (to physicists) thanks
to some explicit constructions using the discrete series representation
of the $N{=}(2,2)$ superconformal algebra \cite{\gp}.  The results
are quite illuminating:  on the one hand, the locally defined map
$${\Cal K}_{\C}(X)\dasharrow {\Cal M}_{\text{complex}}(Y)
\tag{\mirrormap}$$
does {\it not}\/ in general extend throughout ${\Cal K}_{\C}(X)$, but
instead there are points where the theories become singular, and the
map encounters difficulties beyond those points.\footnote{This
phenomenon is already visible in the example
considered in \cite{\CdGP}.}
On the other hand, the image of \eqref{\mirrormap} is {\it not}\/ all
of ${\Cal M}_{\text{complex}}(Y)$---as we have already suggested,
${\Cal K}_{\C}(X)$ is much smaller than ${\Cal M}_{\text{complex}}(Y)$.
This means that there must be
a way to analytically continue the conformal field theories on $X$
beyond the theories specified by ${\Cal K}_{\C}(X)$ (since such theories
occur in ${\Cal M}_{\text{complex}}(Y)$).  This second conclusion was
independently reached by Witten \cite{\phases} on somewhat different
grounds.

What, then, lies beyond the K\"ahler cone for such theories?  In some
cases, the conformal field theories are $\sigma$-models on other
Calabi--Yau manifolds which are obtained by a simple topological surgery
from $X$ (see \cite{\AGM, \small} and \cite{\phases},
or for a more mathematical account,
\cite{\beyond}).
In these cases, as the K\"ahler class is varied and allowed to approach
a wall of the K\"ahler cone, a finite number of holomorphic $2$-spheres
have their areas approach $0$.  When the K\"ahler class is pushed
beyond that wall, the areas of those $2$-spheres would apparently become
negative.  However, the analytically continued $\sigma$-model should instead
 be formulated as
 a $\sigma$-model on a modified manifold $X'$, which is obtained from
$X$ by a surgery along the $2$-spheres
 in such a way that the sign of their (common) homology class has
been reversed (cf.\ \cite{\GS}).

The collection of complexified K\"ahler cones of the various topological
models produces a rich combinatorial structure of regions in the
moduli space corresponding to the different models.  But even these
do not fill up the entire conformal field theory moduli space---there
are additional regions whose associated conformal field theories must
be described by constructions other than
$\sigma$-models \cite{\phases, \small}.
These theories are currently under active study.

\head References \endhead

\def\bysame{\leavevmode\hbox to3em{\hrulefill}\thinspace}

\vskip.25cm

\item{[\AGM]}
P. S. Aspinwall, B. R. Greene, and D. R. Morrison,
{\it {C}alabi--{Y}au moduli space,
  mirror manifolds and spacetime topology change in string theory}, Nuclear
  Phys. B {\bf 416} (1994), 414--480.

\item{[\small]}
\bysame,
{\it Measuring small distances in
  {$N{=}2$} sigma models},
 Nuclear Phys. B {\bf 420} (1994), 184--242.

\item{[\AMone]}
P. S. Aspinwall and D. R. Morrison,
{\it Topological field theory and rational curves},
  Comm. Math. Phys. {\bf 151} (1993), 245--262.

\item{[\AMtwo]}
\bysame,
{\it Chiral rings do not suffice: {$N{=}(2,2)$}
  theories with nonzero fundamental group}, Phys. Lett. B {\bf 334} (1994),
  79--86.

\item{[\AMthree]}
\bysame,
{\it String theory on {K3} surfaces}, Essays on Mirror Manifolds II,
(B. R. Greene and S.-T. Yau, eds.), International Press, Hong Kong,
to appear.

\item{[\BE]}
M. Berger and D. Ebin, {\it Some decompositions of the space of symmetric
tensors on a Riemannian manifold}, J. Differential Geom. {\bf 3} (1969),
 379--392.

\item{[\Besse]}
A. L. Besse, Einstein Manifolds, Springer-Verlag, Berlin, Heidelberg, 1987.

\item{[\Bogomolov]}
F.~A. Bogomolov,
{\it {H}amiltonian {K}\"ahler manifolds},
 Dokl. Akad. Nauk SSSR {\bf 243}, no.  5 (1978), 1101--1104.

\item{[\Calabi]}
 E.~Calabi,
{\it  On {K\"a}hler manifolds with vanishing canonical class},
Algebraic Geometry and Topology, A
  Symposium in Honor of {S}. {L}efschetz (R.~H. Fox et~al., eds.),
 Princeton  University Press, Princeton, 1957, pp. 78--89.

\item{[\CdGP]}
P.~Candelas, X.~C. de~la Ossa, P.~S. Green, and L.~Parkes, {\it A pair of
  {C}alabi--{Y}au manifolds as an exactly soluble superconformal theory},
  Nuclear Phys. B {\bf 359} (1991), 21--74.

\item{[\cls]}
P.~Candelas, M.~Lynker, and R.~Schimmrigk, {\it {C}alabi--{Y}au manifolds in
  weight\-ed ${\Bbb P}_4$}, Nuclear Phys. B {\bf 341} (1990), 383--402.

\item{[\DS]}
M. Dine and N. Seiberg, {\it Microscopic knowledge from macroscopic
physics in string theory}, Nuclear Phys. B {\bf 301} (1988), 357--380.

\item{[\dixon]}
L.~J. Dixon, {\it Some world-sheet properties of superstring compactifications,
  on orbifolds and otherwise}, Superstrings, Unified Theories, and Cosmology
  1987  (G.~Furlan et~al., eds.), World
  Scientific,
1988, pp.~67--126.

\item{[\voa]}
I. Frenkel, J. Lepowsky, and A. Meurman, Vertex Operator Algebras
and the Monster, Academic Press, 1988.

\item{[\GSW]}
M. B. Green, J. H. Schwarz, and E. Witten, Superstring Theory, 2 vols.,
Cambridge University Press, 1987.

\item{[\gp]}
B.~R. Greene and M.~R. Plesser, {\it Duality in {C}alabi--{Y}au moduli space},
  Nuclear Phys. B {\bf 338} (1990), 15--37.

\item{[\Gromov]}
M. Gromov, {\it Soft and hard symplectic geometry},
Proc. Intern. Congress Math.
Berkeley 1986, vol.~1,
American Mathematical Society,
1987, pp.~81--98.

\item{[\GS]}
 V.~Guillemin and S.~Sternberg,
{\it Birational equivalence in the symplectic category},
 Invent. Math. {\bf 97} (1989),  485--522.

\item{[\lvw]}
W.~Lerche, C.~Vafa, and N.~P. Warner, {\it Chiral rings in {$N=2$}
  superconformal theories}, Nuclear Phys. B {\bf 324} (1989), 427--474.

\item{[\MS]}
 D.~McDuff and D.~Salamon,
 $J$-holomorphic Curves and Quantum Cohomology, University Lecture Series,
vol. 6, American Mathematical Society,   1994.

\item{[\guide]}D. R. Morrison,
{\it Mirror symmetry and rational curves on quintic threefolds: A guide for
  mathematicians}, J. Amer. Math. Soc. {\bf 6} (1993), 223--247.

\item{[\compact]}
\bysame,
{\it Compactifications of moduli spaces inspired by mirror symmetry},
  Journ\'ees de G\'eom\'etrie Alg\'ebrique d'Orsay (Juillet 1992),
  Ast\'erisque, vol. 218, Soci\'et\'e Math\-\'ema\-tique de France, 1993,
  pp.~243--271.

\item{[\beyond]}
\bysame,
{\it Beyond the {K\"a}hler cone},
Proc. Hirzebruch's 65th Birthday Workshop
in Algebraic Geometry, to appear.

\item{[\predictions]}
\bysame,
{\it Making enumerative predictions by means of mirror symmetry},
Essays on Mirror Manifolds II,
(B. R. Greene and S.-T. Yau, eds.), International Press, Hong Kong,
to appear.

\item{[\Narain]}
K.~S. Narain,
 {\it New heterotic string theories in uncompactified dimensions ${}<10$},
 Phys. Lett. B {\bf 169} (1986), 41--46.

\item{[\NSW]}
K.~S. Narain, M.~H. Sarmadi, and E.~Witten,
 {\it A note on toroidal compactification of heterotic string
  theory},
 Nuclear Phys. B {\bf 279} (1987), 369--379.

\item{[\RT]}
 Y.~Ruan and G.~Tian,
{\it A mathematical theory of quantum cohomology}, Math. Res. Lett.
{\bf 1}  (1994),  269--278.

\item{[\Tian]}
 G.~Tian,
{\it Smoothness of the universal deformation space of compact {C}alabi--{Y}au
  manifolds and its {P}eterson--{W}eil metric},
Mathematical Aspects of String Theory (S.-T. Yau, ed.),
World Scientific, Singapore,  1987, pp.  629--646.

\item{[\Todorov]}
 A.~N. Todorov,
{\it The {W}eil-{P}etersson geometry of the moduli space of {$SU(n{\ge}3)$}
  ({C}alabi--{Y}au) manifolds, {I}},
 Comm. Math. Phys. {\bf 126} (1989), 325--246.

\item{[\Viehweg]}
E.~Viehweg,
{\it Weak positivity and the stability of certain {H}ilbert points, {III}},
 Invent. Math. {\bf 101} (1990), 521--543.

\item{[\WBerkeley]}
E. Witten, {\it Geometry and physics},
Proc. Intern. Congress Math.
Berkeley 1986, vol.~1,
American Mathematical Society,
1987, pp.~267--303.

\item{[\topsigma]}
\bysame,
{\it Topological sigma models}, Comm.
Math. Phys. {\bf 118} (1988), 411--449.

\item{[\topgrav]}
\bysame,
{\it On the structure of the topological phase of two-dimensional gravity},
 Nuclear Phys. B {\bf 340} (1990),  281--332.

\item{[\Wmirror]}
\bysame, {\it Mirror manifolds and topological field theory},
Essays on Mirror Manifolds,
(S.-T. Yau, ed.), International Press, Hong Kong, 1992, pp. 120--159.

\item{[\phases]}
\bysame,
{\it Phases of {$N{=}2$} theories in two dimensions},
 Nuclear Phys. B {\bf 403} (1993),  159--222.

\item{[\YauICM]}
S.-T. Yau,
{\it The role of partial differential equations in differential geometry},
Proc. Intern. Congress Math.
Helsinki 1978, vol.~1,
Academia Scientiarum Fennica,
1980, pp.~237--250.

\enddocument